\documentclass[11pt]{amsart}
\usepackage{stmaryrd}
\usepackage{xcolor}
\usepackage{epsfig}
\usepackage{enumerate}
\usepackage{amsmath}
\usepackage{amssymb}
\usepackage{amscd}
\usepackage{graphicx}
\usepackage{lineno}
\usepackage{pstricks}
\usepackage{tocvsec2}
\usepackage{mathtools}
\usepackage[noadjust]{cite}
\allowdisplaybreaks[4]

\usepackage [backref, colorlinks, citecolor=red, anchorcolor=black, linkcolor=blue]{hyperref}

\usepackage{amsfonts}
\usepackage[top=35mm, bottom=35mm, left=30mm, right=30mm]{geometry}
\usepackage{verbatim}
\usepackage{graphicx, amssymb, amsmath,amsthm}
\allowdisplaybreaks[4]
\usepackage{appendix}
\usepackage{mathrsfs}
\usepackage{cite}

\makeatletter
\@namedef{subjclassname@2020}{%
	\textup{2020} Mathematics Subject Classification}
\makeatother

\usepackage{xcolor}

\makeatletter
\NewCommandCopy\@@pmod\pmod
\DeclareRobustCommand{\pmod}{\@ifstar\@pmods\@@pmod}
\def\@pmods#1{\mkern4mu({\operator@font mod}\mkern 6mu#1)}
\makeatother

\usepackage{xcolor}

\newtheorem{lemma}{Lemma}[section]

\newtheorem{theorem}{Theorem}[section]

\numberwithin{equation}{section}

\numberwithin{equation}{section}
\renewcommand*{\backref}[1]{}

\makeatletter

\def\MRbibitem{\@ifnextchar[\my@lbibitem\my@bibitem}

\def\mybiblabel#1#2{\@biblabel{{\hyperref{http://www.ams.org/mathscinet-getitem?mr=#1}{}{}{#2}}}}

\def\myhyperanchor#1{\Hy@raisedlink{\hyper@anchorstart{cite.#1}\hyper@anchorend}}

\def\my@lbibitem[#1]#2#3#4\par{%
	\item[\mybiblabel{#2}{#1}\myhyperanchor{#3}\hfill]#4%
	\@ifundefined{ifbackrefparscan}{}{\BR@backref{#3}}%
	\if@filesw{\let\protect\noexpand\immediate
		\write\@auxout{\string\bibcite{#3}{#1}}}\fi\ignorespaces%
}

\def\my@bibitem#1#2#3\par{%
	\refstepcounter\@listctr
	\item[\mybiblabel{#1}{\the\value\@listctr}\myhyperanchor{#2}\hfill]#3%
	\@ifundefined{ifbackrefparscan}{}{\BR@backref{#2}}%
	\if@filesw\immediate\write\@auxout
	{\string\bibcite{#2}{\the\value\@listctr}}\fi\ignorespaces%
}

\makeatother

\newtheorem{conjecture}{Conjecture}[section]


\author{Zhi-Wei Sun}
\address[Zhi-Wei Sun]
{School of Mathematics, Nanjing University, Nanjing, Jiangsu, 210093, P.R. China}
\email{zwsun@nju.edu.cn}

\author{Lilu Zhao}
\address[Lilu Zhao]
{School of Mathematical Sciences, University of Science and Technology of China, Hefei, Anhui, 230026, P.R. China}
\email{zhaolilu@ustc.edu.cn}

\begin{document}

	\begin{abstract}
		Let $G$ be a finite group with $|G|=p^m$ where $p$ is a prime and $m$ is a positive integer. Let $k<p$. Let $a_1,\ldots,a_k\in G$ be pairwise distinct and let $b_1,\ldots,b_k\in G$. Then there exists a permutation $\sigma$ on $1,\ldots,k$ such that $a_1b_{\sigma(1)},\ldots,a_kb_{\sigma(k)}$ are pairwise distinct. This extends a theorem of Feng, Sun and Xiang, who proved that the conclusion holds in abelian $p$-groups.
	\end{abstract}

	
	\title{Exterior Algebra and an Extension of the Feng-Sun-Xiang Theorem in $p$-groups}
	\maketitle

	\parskip 10pt         

	
	\section{Introduction}

Let $G=\{a_1,\ldots,a_n\}$ be an abelian group, and let $b_1,\ldots,b_n\in G$. In 1952, M. Hall, Jr. \cite{Hall} proved that there exists a permutation $\sigma\in S_n$ such that $a_1b_{\sigma(1)},\ldots,a_kb_{\sigma(k)}$ are pairwise distinct if and only if $b_1\cdots b_n=e$, where $e$ denotes the identity element of $G$. Snevily \cite{Snevily} proposed an interesting conjecture that if $A=\{a_1,\ldots,a_k\}$ and $B=\{b_1,\ldots,b_k\}$ are subsets of an abelian group $G$ of odd order then $a_1b_{\sigma(1)},\ldots,a_kb_{\sigma(k)}$ are pairwise distinct for some permutation $\sigma\in S_k$. Snevily's conjecture was confirmed by Alon \cite{Alon} when $G$ is a cyclic group of prime order. Later, Dasgupta, K\'arolyi, Serra and Szegedy \cite{DKSS} confirmed Snevily's conjecture when $G$ is a cyclic group of odd order. Gao and Wang \cite{GW} verified Snevily's conjecture under the condition that $k<\sqrt{p(G)}$, where $p(G)$ denotes the smallest prime factor of $|G|$. Snevily's conjecture was finally solved by Arsovski \cite{Arsovski}. For the development of Snevily's conjecture in combination with restricted sumsets, one may refer to \cite{Sun}.

When $G$ is a cyclic group of prime order, Alon \cite{Alon} actually proved a stronger result, in which $b_1,\ldots,b_k$ are not necessarily distinct.

	\begin{theorem}[Alon \cite{Alon}]\label{alon}Let $G$ be a cyclic group of prime order $p$. Let $k<p$. Let $a_1,\ldots,a_k\in G$ be pairwise distinct and let $b_1,\ldots,b_k\in G$. Then there exists a permutation $\sigma$ on $1,\ldots,k$ such that $a_1b_{\sigma(1)},\ldots,a_kb_{\sigma(k)}$ are pairwise distinct.
	\end{theorem}
 	Throughout this paper, we assume $p$ is a prime number. Dasgupta, K\'arolyi, Serra and Szegedy \cite{DKSS} extended Alon’s result (Theorem \ref{alon}) to cyclic groups of prime power order as well as elementary abelian groups.
 	\begin{theorem}[Dasgupta-K\'arolyi-Serra-Szegedy \cite{DKSS}]\label{dkss}Let $G$ be a cyclic group of order $p^{m}$ with $m\in \mathbb{Z}^{+}$ or the elementary abelian p-group. Let $k<p$. Let $a_1,\ldots,a_k\in G$ be pairwise distinct and let $b_1,\ldots,b_k\in G$. Then $a_1b_{\sigma(1)},\ldots,a_kb_{\sigma(k)}$ are pairwise distinct for some permutation $\sigma\in S_k$.
	\end{theorem}
 
 Motivated by Theorems \ref{alon} and \ref{dkss}, the following conjecture was proposed in \cite{DKSS}.

  	\begin{conjecture}[Dasgupta-K\'arolyi-Serra-Szegedy]\label{dkssconj}Let $G$ be a finite abelian group with $|G|>1$. Let $p$ be the smallest prime of $|G|$. Let $k<p$. Let $a_1,\ldots,a_k\in G$ be pairwise distinct and let $b_1,\ldots,b_k\in G$. Then $a_1b_{\sigma(1)},\ldots,a_kb_{\sigma(k)}$ are pairwise distinct for some permutation $\sigma\in S_k$.
	\end{conjecture}

In 2011, Feng, Sun and Xiang \cite{FSX} introduced the method of exterior algebra to study Conjecture \ref{dkssconj} and related problems. They proved the following (see \cite[Corollary 1.7]{FSX}).

	\begin{theorem}[Feng-Sun-Xiang \cite{FSX}]\label{FSX}Conjecture \ref{dkssconj} is true when $G$ is an abelian $p$-group.
	\end{theorem}
Indeed, Feng, Sun and Xiang obtained a stronger conclusion which asserts that Conjecture \ref{dkssconj} is true when the
smallest prime divisor of $|G|$ is greater than $k$ and any other prime divisor of $G$ (if there is
any) is greater than $k!$. Feng, Sun and Xiang made use of the character theory of abelian groups, which is not applicable to non-abelian groups.

The purpose of this paper is to further develop the method of exterior algebra into this topic, and as an application we extend Theorem \ref{FSX} to finite (non-abelian) $p$-groups. The main result in this paper is the following.
	\begin{theorem}\label{main}Let $G$ be a finite $p$-group. Let $k<p$. Let $a_1,\ldots,a_k\in G$ be pairwise distinct and let $b_1,\ldots,b_k\in G$. Then there exists a permutation $\sigma$ on $1,\ldots,k$ such that $a_1b_{\sigma(1)},\ldots,a_kb_{\sigma(k)}$ are pairwise distinct.
	\end{theorem}
	
	In the next section we will give some background on exterior algebras. Theorem \ref{main} will be proved in Section 3.
	
	\section{Preparations and exterior algebras}
	
	Let $\mathbb{F}$ be a field and let $V$ be a vector space over $\mathbb{F}$. Let $k\geqslant 2$ be a positive integer. Let $V^k$ denote the $k$-fold direct sum of $V$. The $k$-th exterior power of $V$, denoted by $\bigwedge^{k}V$, is with elements in the form
$$v_1\wedge \cdots \wedge v_k,$$
where $v_1,\ldots,v_k\in V$.
For an alternating multilinear map $\widetilde{h}:\ V^k\rightarrow \bigwedge^{k}V$, there exists a unique linear map $h:\ \bigwedge^{k}V\rightarrow \bigwedge^{k}V$ satisfying
$$h(v_1\wedge \cdots \wedge v_k)=\widetilde{h}(v_1,\ldots,v_k).$$
Throughout this section, we fix 
$$w_1,\ldots,w_k\in V.$$
Now we introduce $\widetilde{h}:\ V^k\rightarrow \bigwedge^{k}V$ (depending on $w_1,\ldots,w_k$) defined by
\begin{align}\label{definehtilde}\widetilde{h}(v_1,\ldots,v_k)=\sum_{\sigma\in S_k}v_1w_{\sigma(1)} \wedge \cdots \wedge v_kw_{\sigma(k)}.\end{align}
This is a multilinear map. It is worth pointing out that the map
$$(v_1,\ldots,v_k)\rightarrow v_1w_1 \wedge \cdots \wedge v_kw_k$$
is also multilinear, but not alternating. However, $\widetilde{h}$ is an alternating multilinear map.
\begin{lemma}\label{lemma1}Let $\widetilde{h}:\ V^k\rightarrow \bigwedge^{k}V$ be defined in \eqref{definehtilde}. Then $\widetilde{h}$ is an alternating multilinear map.
In particular, it induces a linear map $h:\ \bigwedge^{k}V\rightarrow \bigwedge^{k}V$ satisfying
\begin{align}\label{defineh}h(v_1 \wedge \cdots \wedge v_k)=\sum_{\sigma\in S_k}v_1w_{\sigma(1)} \wedge \cdots \wedge v_kw_{\sigma(k)}.\end{align}
  \end{lemma}	
\begin{proof}It is easy to see that $\widetilde{h}$ is multilinear. It suffices to prove the alternating property. Suppose that $1\leqslant i<j\leqslant k$. We need to prove that if $v_i=v_j$ then
\begin{align}\label{alternating}\widetilde{h}(v_1,\ldots,v_i,\ldots,v_j,\ldots, v_k)=0.\end{align}

We introduce the transposition $\tau=(i,j)\in S_k$. For any $\sigma\in S_k$, noting that $\sigma \tau^2=\sigma$ and $\sigma\tau\not=\sigma$, we can view $\{\sigma,\sigma \tau\}$ as a pair. Then $S_k$ can be divided into $\frac{k!}{2}$ such pairs. Note that 
$$\sigma\tau(j)=\sigma(i),\ \sigma\tau(i)=\sigma(j)$$
and 
$$\sigma\tau(\ell)=\sigma(\ell)\ \ \textrm{for any }\ \ell\not\in \{i,j\}.$$
Then we have
\begin{align*}v_1w_{\sigma\tau(1)} \wedge \cdots \wedge v_kw_{\sigma\tau(k)}
= &\, v_1w_{\sigma(1)} \wedge \cdots \wedge v_iw_{\sigma(j)}\wedge \cdots \wedge v_jw_{\sigma (i)} \wedge \cdots \wedge v_kw_{\sigma(k)}
\\= &\, -v_1w_{\sigma(1)} \wedge \cdots \wedge v_jw_{\sigma(i)}\wedge \cdots \wedge v_iw_{\sigma (j)} \wedge \cdots \wedge v_kw_{\sigma(k)}.\end{align*}
Since $v_i=v_j$, we further have 
\begin{align*}v_1w_{\sigma\tau(1)} \wedge \cdots \wedge v_kw_{\sigma\tau(k)}=-v_1w_{\sigma(1)} \wedge \cdots \wedge v_iw_{\sigma(i)}\wedge \cdots \wedge v_jw_{\sigma (j)} \wedge \cdots \wedge v_kw_{\sigma(k)}.\end{align*}
Thus, 
for each pair $\{\sigma,\sigma \tau\}$, we have
\begin{align*}v_1w_{\sigma(1)} \wedge \cdots \wedge v_kw_{\sigma(k)}+v_1w_{\sigma\tau(1)} \wedge \cdots \wedge v_kw_{\sigma\tau(k)}=0.\end{align*}
Summing over all pairs, we establish \eqref{alternating}. By the universal property, $h:\ \bigwedge^{k}V\rightarrow \bigwedge^{k}V$  in \eqref{defineh} is a linear map.
\end{proof}	

For a prime $p$, we use $\mathbb{F}_p$ to denote the prime field. Let $G$ be a finite group. We shall focus on the group algebra
\begin{align}\label{FG}V=\mathbb{F}_p[G].\end{align}
The elements in $V$ will be denoted by
$$\sum_{g\in G}c_g g$$
with $c_g\in \mathbb{F}_p$. Let $G=\{g_1,\ldots,g_{|G|}\}$. The dimension of $V$ is $\binom{|G|}{k}$ with a basis consisting of
\begin{align}\label{basis}g_{j_1} \wedge \cdots \wedge g_{j_k},\end{align}
where $1\leqslant j_1< \cdots<j_k\leqslant |G|$.
The augmentation ideal of $\mathbb{F}_p[G]$ is
\begin{align}\label{ideal}I=\Big\{\sum_{g\in G}c_g g:\ \sum_{g\in G}c_g =0\Big\}.\end{align}
The following is a classical result.
\begin{lemma}\label{lemma2} Let $G$ be a finite $p$-group, that is $|G|=p^m$ for some positive integer $m$. Let $I$ be the augmentation ideal of $\mathbb{F}_p[G]$. Then
$I$ is nilpotent, that is $I^{(n)}=0$ for some positive integer $n$. In particular, for any $u_1,\ldots,u_n\in I$, one has $u_1\cdots u_n=0$.
\end{lemma} 
\begin{proof}See \cite[Lemma 1.6]{Passman}.\end{proof}

	\section{Proof of Theorem \ref{main}}

\noindent {\it Proof of Theorem \ref{main}.} We choose $V$ as in \eqref{FG}. We introduce $f:\ \bigwedge^{k}V\rightarrow \bigwedge^{k}V$ defined by
\begin{align}\label{definef}f(v_1 \wedge \cdots \wedge v_k)=\sum_{\sigma\in S_k}v_1b_{\sigma(1)} \wedge \cdots \wedge v_kb_{\sigma(k)}.\end{align}
By Lemma \ref{lemma1}, this is a well-defined linear map. Our aim is to prove that $f$ is invertible. 

Let $\mathcal{J}=\{0,1\}^{k}\setminus \{\mathbf{0}\}\subseteq \mathbb{Z}^k$, i.e. $\mathcal{J}$ is the set of all nonzero $k$-dimensional vectors with all coordinates equal to $0$ or $1$. We can express $v_1b_{\sigma(1)} \wedge \cdots \wedge v_kb_{\sigma(k)}$ in the form
\begin{align*}v_1b_{\sigma(1)} \wedge \cdots \wedge v_kb_{\sigma(k)}=&\, v_1(b_{\sigma(1)}-e+e) \wedge \cdots \wedge v_k(b_{\sigma(k)}-e+e)
\\= &\, \sum_{(j_1,\ldots,j_k)\in \mathcal{J}}v_1(b_{\sigma(1)}-e)^{j_1} \wedge \cdots \wedge v_k(b_{\sigma(k)}-e)^{j_k}+v_1 \wedge \cdots \wedge v_k,\end{align*}
where $e$ is the identity element of the group $G$.
Then from \eqref{definef}, we have
\begin{align*}f(v_1 \wedge \cdots \wedge v_k)= \sum_{\sigma\in S_k}\sum_{(j_1,\ldots,j_k)\in \mathcal{J}}v_1(b_{\sigma(1)}-e)^{j_1} \wedge \cdots \wedge v_k(b_{\sigma(k)}-e)^{j_k}+k!v_1 \wedge \cdots \wedge v_k.\end{align*}
We use $\textrm{id}: \ \bigwedge^{k}V\rightarrow \bigwedge^{k}V$ to denote identity map. Let 
$$\delta=f-k!\textrm{id}.$$
 Clearly
\begin{align*}\delta(v_1 \wedge \cdots \wedge v_k)=\sum_{\sigma\in S_k}\sum_{(j_1,\ldots,j_k)\in \mathcal{J}}v_1(b_{\sigma(1)}-e)^{j_1} \wedge \cdots \wedge v_k(b_{\sigma(k)}-e)^{j_k}.\end{align*}
For a vector $\mathbf{j}^{(s)}\in \mathcal{J}$, we adopt the notation $\mathbf{j}^{(s)}=(j_1^{(s)},\ldots, j_k^{(s)})$. Now we rewrite the above in the form
\begin{align*}\delta(v_1 \wedge \cdots \wedge v_k)=\sum_{\sigma_1\in S_k}\sum_{\mathbf{j}^{(1)}\in \mathcal{J}}v_1(b_{\sigma(1)}-e)^{j_1^{(1)}} \wedge \cdots \wedge v_k(b_{\sigma(k)}-e)^{j_k^{(1)}}.\end{align*}

Let 
$$\delta^{r}=\underbrace{\delta\circ\cdots \circ \delta}_{r\ \textrm{times}}.$$
Note that $\delta:\ \bigwedge^{k}V\rightarrow \bigwedge^{k}V$ is a linear map. Then we have
\begin{align*}&\, \delta^r(v_1 \wedge \cdots \wedge v_k)
\\=&\, \sum_{\sigma_1,\ldots,\sigma_{r}\in S_k}\sum_{\mathbf{j}^{(1)},\ldots,\mathbf{j}^{(r)}\in \mathcal{J}}
v_1(b_{\sigma_1(1)}-e)^{j_1^{(1)}}\cdots(b_{\sigma_r(k)}-e)^{j_1^{(r)}} \wedge \cdots \wedge 
v_k(b_{\sigma_1(k)}-e)^{j_k^{(1)}}\cdots (b_{\sigma_r(k)}-e)^{j_k^{(r)}}.\end{align*}

Note that for each $s~(1\leqslant s\leqslant r)$, $\mathbf{j}^{(s)}$ in (a nonzero vector) in $\mathcal{J}$. Therefore, $\sum_{i=1}^kj_i^{(s)}\geqslant 1$ and 
\begin{align}\label{ineq}\sum_{i=1}^k\sum_{s=1}^rj_i^{(s)}=\sum_{s=1}^r\sum_{i=1}^kj_i^{(s)}\geqslant r.\end{align}
Let $I$ be the augmentation ideal of $\mathbb{F}_p[G]$. Then by Lemma \ref{lemma2}, $I^{(n)}=0$ for some positive integer $n$. 
We choose $r=kn$. Then by \eqref{ineq}, there exists $i(1\leqslant i\leqslant k)$ such that 
$$\sum_{s=1}^rj_i^{(s)}\geqslant n.$$
Noting that $b_\ell-e\in I$ for any $1\leqslant \ell \leqslant k$, now we can conclude that
$$(b_{\sigma_1(i)}-e)^{j_i^{(1)}}\cdots (b_{\sigma_r(i)}-e)^{j_i^{(r)}}=0,$$
	and this implies $\delta^r=0$ (a zero map). 

We write $T=k!\textrm{id}$ and $S=-\delta$. Then $f=T-S$. Note that 
\begin{align*}f(T^{r-1}+T^{r-2}S+\cdots+TS^{r-2}+S^{r-1})=&\, (T-S)(T^{r-1}+T^{r-2}S+\cdots+TS^{r-2}+S^{r-1})
\\ = &\, T^r-S^r=T^r=(k!)^r\textrm{id}.\end{align*}
 Now we see that $f$ is invertible provided that $k<p$. 

Since $a_1,\ldots,a_k$ are pairwise distinct elements in $G$, by the argument around \eqref{basis}, we have $a_1 \wedge \cdots \wedge a_k\not=0$. Since $f$ is invertible, we conclude that $f(a_1 \wedge \cdots \wedge a_k)\not=0$.
 In particular, $a_1b_{\sigma(1)} \wedge \cdots \wedge a_kb_{\sigma(k)}\not=0$ for some $\sigma\in S_k$. It follows that  $a_1b_{\sigma(1)},\ldots,a_kb_{\sigma(k)}$ are pairwise distinct.
 
 We complete the proof of Theorem \ref{main}.

\noindent {\bf Acknowledgement}.
	This work is supported by the National Key Research and Development Program of China (No. 2021YFA1000701) and National Natural Science Foundation of China (No. 12371004 and No. 12471088). 
	
	\bibliographystyle{acm}

\end{document}